\documentclass[12pt]{article}

\title{Intrinsically triple-linked graphs in $\mathbb{R}P^3$ \footnote{This
material is based upon work supported by the National Science Foundation under
Grant No. 0646847, and the National Security Administration under Grant No.
42652.}}

\author{
J. Federman \footnote{Department of Mathematics, SUNY Potsdam, Potsdam, NY
13676}
 \and J. Foisy \footnote{Department of Mathematics, SUNY Potsdam, Potsdam, NY
13676}
 \and K. McNamara \footnote{Department of Mathematics \& Statistics, James
Madison University, Harrisonburg, VA 22807}
 \and E.Stark \footnote{Department of Mathematics, Pomona College, Claremont, CA
91711}}

\usepackage{amsmath,amssymb,amsthm}
\usepackage{graphicx,epsfig,color}

\theoremstyle{definition}
\newtheorem{thm}{Theorem}

\newtheorem{lemma}[thm]{Lemma}
\newtheorem{prop}[thm]{Proposition}

\newcommand{\Z}{\mathbb{Z}}

\newcommand{\R}{\mathbb{R}}

\begin{document}

\maketitle

\begin{abstract}
Flapan--Naimi--Pommersheim \cite{FNP} showed that every spatial embedding of $K_{10}$, the complete graph on ten vertices, contains a non-split three-component link; that is, $K_{10}$ is {\it intrinsically triple-linked} in $\R^3$.  The work of Bowlin--Foisy \cite{BF} and Flapan--Foisy--Naimi--Pommersheim \cite{FFNP} extended the list of known intrinsically triple-linked graphs in $\mathbb{R}^3$ to include several other families of graphs.  In this paper, we will show that while some of these graphs can be embedded 3-linklessly in $\mathbb{R}P^3$, $K_{10}$ is intrinsically triple-linked in $\mathbb{R}P^3$.  
\end{abstract}


\section{Introduction}

There is a classic theory of knots and links in Euclidean $3$-space (or the $3$-sphere), and, as Manturov \cite{Man} points out in his book, there is a sympathetic theory of knots and links in $\R P^3$. Drobotukhina \cite{Dr} developed an analog of the Jones polynomial for the case of oriented links in $\R P^3$, and Mroczkowski \cite{M} described a method to unknot knots and links in $\R P^3$ through an analog of classical knot and link diagrams for knots in $\R^3$. Flapan--Howards--Lawrence--Mellor \cite{FHLM} investigate intrinsic linking and knotting in arbitrary $3$-manifolds. Here, following Bustamente et al. \cite{REU07}, we use a weaker notion of unlink than was used in \cite{FHLM}, and we examine the intrinsic linking properties of graphs embedded in $\R P^3$. In particular, we will examine graphs that contain a $3$-component non-split link in every embedding into $\R P^3$.

\vspace{\baselineskip}


Real projective $3$-space $\R P^3$ can be obtained from the $3$-ball $D^3$ by
identifying opposite points of its boundary; hence, a link in
$\R P^3$ consists of a union of arcs and loops so that the endpoints of any arc
lie on antipodal boundary points of the $3$-ball. We may use ambient isotopy to
move all arcs so that their endpoints lie on a fixed great circle, the
``equator'' of the ball. Therefore, a link may be represented in $\R P^2$ by its
projection onto a $2$-disk, $D^2$, whose boundary is the equator, with
antipodal points identified. 

\vspace{\baselineskip}

Projective space has a non-trivial first homology group, $H_1(\R P^3) \cong
\Z_2$. Let $g$, the cycle consisting of a line in $D^3$ running between the
north and south poles, be the generator of this group. Using crossing changes
and ambient isotopy on an $\R P^2$ projection of
a knot, Mroczkowski \cite{M} showed that every knot in $\R P^3$ can be
transformed into
either the trivial cycle or $g$. Thus, there are two non-equivalent unknots in
$\mathbb{R}P^3$.  Cycles that can be unknotted into a cycle homologous to $g$
will be referred to as {\it 1-homologous cycles}.  Cycles that can be
unknotted into a trivial cycle will be referred to as {\it 0-homologous
cycles}. 

\vspace{\baselineskip}

Following \cite{REU07}, we say a link in $\R P^3$ is {\it splittable} if one
component can be contained within a $3$-ball embedded in $\R P^3$, while the
other component lies in the complement of the $3$-ball. Otherwise, a link in
$\R P^3$ is said to be {\it non-split}. A non-split link may be formed in one of
three ways in $\R P^3$: by two $0$-homologous cycles, by a $0$-homologous cycle
and a $1$-homologous cycle, and by two $1$-homologous cycles. Moreover, since a
$1$-homologous cycle cannot be contained within a ball embedded in $\R P^3$,
two disjoint $1$-homologous cycles will always form a non-split link. In this
paper, will we refer to non-split linked cycles as {\it linked cycles} and to
an embedded graph as {\it linked} if it contains a non-split link. 

\vspace{\baselineskip}

A graph
$H$ is a {\it minor} of $G$ if $H$ can be obtained from $G$ through a series of
vertex removals, edge removals, or edge contractions. A graph $G$ is said to be
{\it minor-minimal} with respect to property $P$ if $G$ has property $P$, but no
minor of $G$ has property $P$. The complete set of minor-minimal intrinsically
linked graphs in $\R^3$ is given by the Petersen Family graphs, including $K_6$
and the graphs obtained from $K_6$ by $\Delta-Y$ and $Y-\Delta$ exchanges
\cite{CG, RST, S:84}. However, all Petersen Family graphs except $K_{4,4}-edge$
embed
linklessly in $\R P^3$, as shown in \cite{REU07}, a paper which also exhibits $597$
graphs that are minor-minimal intrinsically linked in $\R P^3$. The complete set
of minor-minimal intrinsically linked graphs in $\R P^3$ is finite \cite{RS},
and remains to be found. 

\vspace{\baselineskip}

A {\it non-split triple-link} is a non-split link of three components, which,
in an abuse of language, will be referred to as a {\it triple-link} in this
paper. An embedding of a graph is {\it triple-linked} if it contains a
non-split link of three components, and a graph is {\it intrinsically triple
linked in $X$}, a topological space, if every embedding of the graph into $X$
contains a non-split triple-link. 

\vspace{\baselineskip}

Conway, Gordon, \cite{CG} and Sachs \cite{S:83, S:84} proved that $K_6$
is intrinsically linked in $\R^3$. In contrast, $K_6$ can be linklessly embedded in $\R P^3$ (see Figure \ref{linkless}). In \cite{REU07}, $7$ is shown to be the smallest $n$ for which $K_n$ is intrinsically linked in $\R P^3$. Flapan--Naimi--Pommersheim \cite{FNP} proved $10$ is the smallest $n$ for which $K_n$ is intrinsically triple-linked in $\R^3$. We show, in Section 3, that $10$ is also the smallest $n$ for which $K_n$ is intrinsically triple-linked in $\R P^3$. It remains to be shown whether $K_{10}$ is minor-minimal with respect to triple-linking in $\R P^3$. 

\vspace{\baselineskip}

In Section 4, we show that two intrinsically triple-linked graphs in $\R^3$ can
be
embedded $3$-linklessly in $\R P^3$, and exhibit two other minor-minimal
intrinsically triple-linked graphs in $\R P^3$. A complete set of minor-minimal
intrinsically triple-linked graphs in both $\R^3$ and $\R P^3$ remains to be
found. Such sets are finite due to the result in \cite{RS}.

\section{Definitions and preliminary lemmas}

We begin with some elementary definitions and notation. A {\it graph}, $G =
(V,E)$, is a set of {\it vertices}, $V(G)$, and {\it edges}, $E(G)$, where an
edge is an unordered pair $(v_1, v_2)$ with $v_1, v_2 \in V(G)$. If $G$ is a
graph with $v_1, \ldots, v_n \in V(G)$ and $(v_1,v_2), (v_2, v_3), \ldots,$ 
$(v_{n-1},v_n), (v_n, v_1) \in E(G)$, with $v_i \neq v_j$ for all $i \neq j$,
then the sequences of vertices  $v_1, \ldots, v_n$ and edges $(v_1,v_2), (v_2,
v_3), \ldots, (v_{n-1},v_n), (v_n, v_1)$ is an {\it $n$-cycle} of $G$, denoted
$(v_1, \ldots, v_n)$. In this paper, we also refer to the image of a cycle
under an embedding as an $n$-cycle. 

\vspace{\baselineskip}
\begin{figure}
\begin{center}
\includegraphics[scale=0.4]{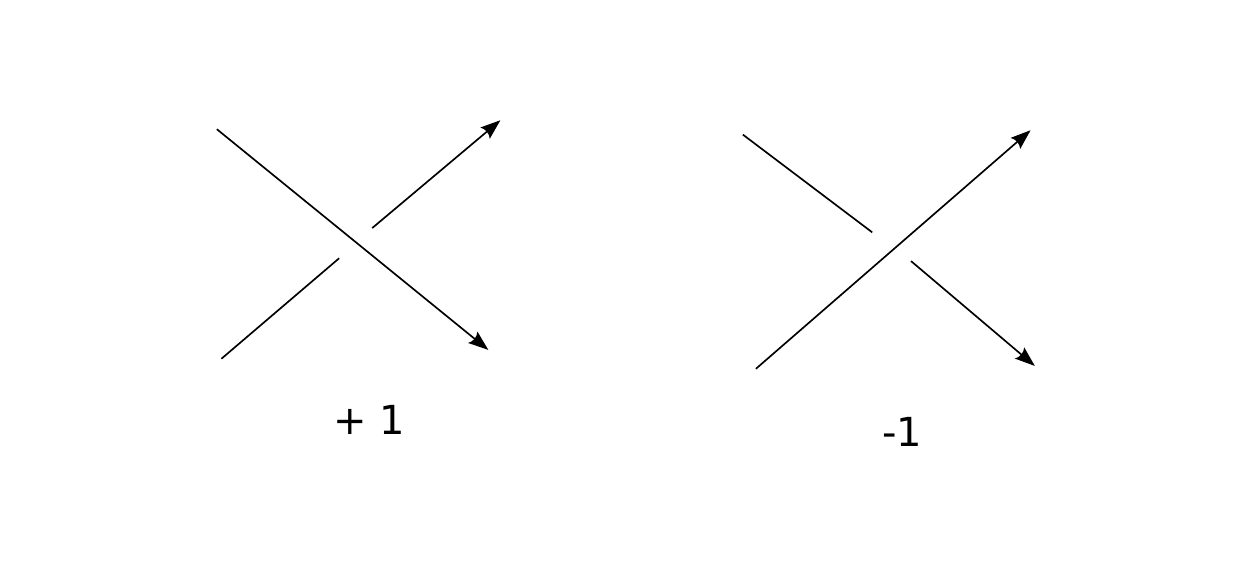}
\vskip -.4 in
\caption{Link crossings.}
\label{linkingnumber}
\end{center}
\end{figure}
\begin{figure}

\begin{center}
\vskip -.5in
\includegraphics[scale=0.3]{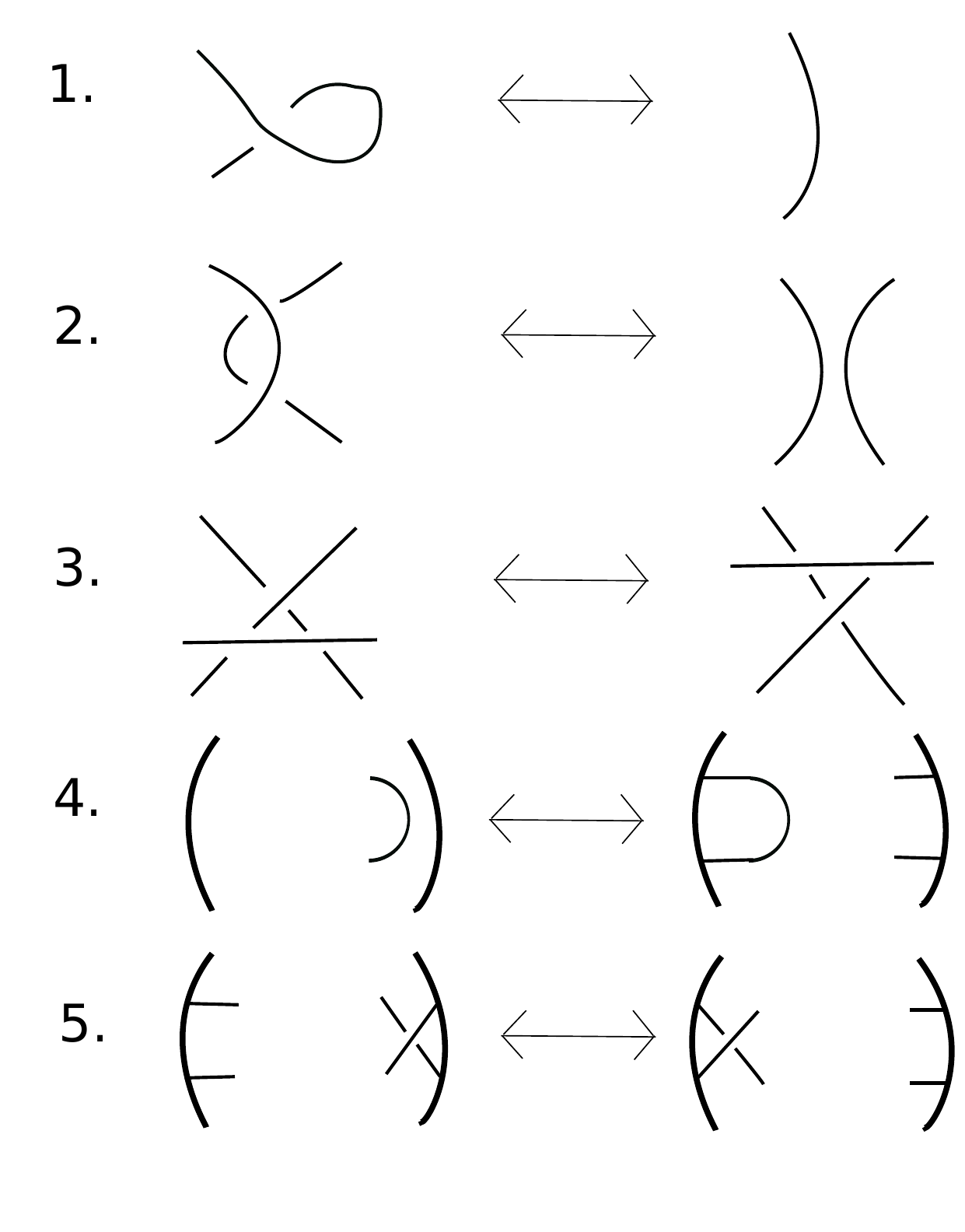}
\vskip -.1 in
\caption{Generalized Reidemeister moves in $\R P^3$.}
\label{Reid}
\end{center}
\end{figure}

If $G$ is a graph and $v_1, \ldots, v_n \in V(G)$, define the {\it induced
subgraph}, $G[v_1, \ldots, v_n]$, to be the subgraph of $G$ with 
\begin{eqnarray*}
V(G[v_1, \ldots, v_n]) &=& \{v_1, \ldots, v_n\}\\
E(G[v_1, \ldots, v_n]) &=& \{(v_i, v_j) \in E(G)| v_i, v_j \in \{v_1, \ldots,
v_n\}\}. 
\end{eqnarray*}

The classical notion of linking number extends to links embedded in $\R P^3$.
Suppose $L$ and $K$ are two loops embedded in $\R P^3$; orient $L$ and $K$. At
each crossing, assign $+1$ or $-1$ as drawn in Figure \ref{linkingnumber}. Then
the {\it mod $2$ linking number of $L$ and $K$}, $lk(L,K)$, is
the sum of the numbers, $+1$ or $-1$, at each crossing in the embedding of $L$
and $K$ divided by two, taken modulo $2$. In $\R P^3$, there are five
generalized Reidemeister moves, which are drawn in Figure \ref{Reid} \cite{Man}.
As in $\R^3$, one can use Reidemeister moves to justify that mod $2$ linking
number is well-defined in $\R P^3$.
In particular, the mod $2$ linking number of a splittable two-component link is
$0$.
However, in $\R P^3$, the mod $2$ linking number need not be an integer; for
example,
two disjoint $1$-homologous cycles can have mod $2$ linking number $\pm
\frac{1}{2}$. 

\vspace{\baselineskip}

In this paper, the following lemmas provide us information about carefully
chosen induced subgraphs of the graphs we study. 

\begin{lemma} \label{K7} \cite{REU07} {\it The graphs obtained by removing two
edges from $K_7$ and removing one edge from $K_{4,4}$ are intrinsically linked
in $\mathbb{R}P^3$.}
\end{lemma}

\begin{lemma} \label{OHK4} \cite{REU07} {\it Given a linkless embedding of $K_6$
in $\mathbb{R}P^3$, no $K_4$ subgraph can have all 0-homologous cycles.}
\end{lemma}

In addition, we use the following elementary observation.

\begin{lemma} \label{evenK4} {\it For every embedding into $\mathbb{R}P^3$,
$K_4$ has an even number of 1-homologous 3-cycles.}
\end{lemma}

The following two lemmas were shown true in $\R^3$ by \cite{FNP} and \cite{BF},
respectively. In each case, the proof holds analogously in $\R P^3$. 

\begin{lemma} \label{LL} {\it Let $G$ be a graph embedded in $\mathbb{R}P^3$
that contains cycles $C_1$, $C_2$, $C_3$ and $C_4$. Suppose $C_1$ and $C_4$ are
disjoint from each other and from $C_2$ and $C_3$ and suppose $C_2 \cap C_3$ is
a simple path. If $lk(C_1,C_2) \neq 0$ and $lk(C_3,C_4) \neq 0$, then $G$
contains a non-split three-component link.}
\end{lemma}

\begin{lemma} \label{BFLL} {\it In an embedded graph with mutually
disjoint
simple closed curves, $C_1$, $C_2$, $C_3$, and $C_4$, and two disjoint paths
$x_1$ and $x_2$ such that $x_1$ and $x_2$ begin in $C_2$ and end in $C_3$, if
$lk(C_1,C_2) \neq 0$  and $lk(C_3,C_4) \neq 0$, then the embedded
graph contains a non-split three component link.}
\end{lemma}

\section{Intrinsically triple-linked complete graphs on $n$ vertices}

The first proposition of this section, that $K_{11}$ is intrinsically
triple-linked in $\R
P^3$, is not the main result of this paper. In fact, our main result, that
$K_{10}$ is intrinsically triple-linked in $\R P^3$, implies this proposition,
by a result of \cite{NT}. However, the proof is included because it is
(relatively) concise and follows from examining four carefully chosen subgraphs
of $K_{11}$ and applying Lemmas \ref{LL} and \ref{BFLL}. 

\begin{prop} {\it The graph $K_{11}$ is intrinsically triple-linked in
$\mathbb{R}P^3$.}
\end{prop}

\begin{proof}
 Let $G$ be a complete graph isomorphic to $K_{11}$ with vertices labeled
$\{1,2, \ldots, 11\}$. Embed $G$ in $\R P^3$. 

\vspace{\baselineskip}

Since $K_7$ is intrinsically linked in $\R P^3$, the graph $G[1,2,3,4,5,6,7]
\cong K_7$ contains a pair of linked cycles. Without loss of generality,
suppose the linked cycles are $C_1 = (1,2,3)$ and $C_2' = (4,5,6,7)$. Homologically, the cycle $(4,5,6,7)$ is the sum of the cycles $(4,5,6)$ and $(4,6,7)$. Thus, $lk((1,2,3), (4,5,6,7)) = lk((1,2,3), (4,5,6)) + lk((1,2,3), (4,6,7))$. Since the numbers on the right-hand side cannot both equal zero, without loss of generality, $C_1 = (1,2,3)$ links with $C_2 = (4,5,6)$.


\vspace{\baselineskip}

Again, since $K_7$ is intrinsically linked in $\R P^3$, the subgraph
$G[5,6,7,8,9,10,11] \cong K_7$ contains a pair of linked cycles. In the manner
described above, this pair of cycles may be reduced to two linked $3$-cycles.
If it is not the case that one cycle contains $\{5\}$ and one cycle contains
$\{6\}$, then Lemma \ref{LL} applies, and $G$ is triple-linked. To handle the
other case, suppose, without loss of generality, that $C_3 = (5,7,9)$ and $C_4
= (6,8,10)$ are the pair of linked cycles in $G[5,6,7,8,9,10,11]$. 

\vspace{\baselineskip}

To obtain two collections of disjoint $1$-homologous cycles, consider two
subgraphs isomorphic to $K_6$. First, if $G[1,2,3,4,6,11] \cong K_6$ contains a
pair of linked cycles, then one cycle shares vertex $\{6\}$ with $C_4$ and both
are disjoint from $C_3$, so Lemma \ref{LL} applies and $G$ is triple-linked.
Otherwise, by Lemma \ref{OHK4}, the set $A = \{(1, 2, 3), (1, 2, 11), (1, 3,
11), (2, 3, 11)\}$ contains a $1$-homologous cycle, $C_5$. 

\vspace{\baselineskip}

Similarly, if $G[6,7,8,9,10,11] \cong K_6$ contains a pair of linked cycles,
then one cycle shares vertex $\{6\}$ with $C_2$ and both are disjoint from
$C_1$. So, Lemma \ref{LL} applies and $G$ is triple-linked. Otherwise, by Lemma
\ref{OHK4}, the set $B = \{ (7, 8, 9), (7, 8, 10), (7, 9, 10)$, $(8, 9, 10)\}$
contains a $1$-homologous cycle, $C_6$. 

\vspace{\baselineskip}

Since $A \cap B = \emptyset$, $C_5 \in A$ and $C_6 \in B$ are disjoint
$1$-homologous
cycles and hence linked. So, $C_2$ and $C_6$ are disjoint from
each other and from $C_1$ and $C_5$. In the case that $C_1 = C_5$, $C_1,
C_2$, and $C_6$ form a triple link. Otherwise, $C_1 \cap C_5$ is a simple
path, so $G$ contains a triple-link by Lemma \ref{LL}.
\end{proof}

To prove that $K_{10}$ is intrinsically triple-linked in $\R P^3$, we first
describe how its subgraphs isomorphic to $K_6$ must be embedded. 

\begin{prop} \label{TD}
 If $G$ is isomorphic to $K_6$ and embedded in $\R P^3$ and $G$ contains two
disjoint $0$-homologous cycles, then $G$ contains a non-split link.
\end{prop}
\begin{proof}
 Let $G$ be isomorphic to $K_6$ and suppose $G$ is embedded so that it
contains two disjoint $0$-homologous cycles and no non-split link. Without loss
of generality, let $(1,2,3)$ and $(4,5,6)$ be $0$-homologous cycles in $G$.
Consider $G[1,2,3,4]$. Since $G$ is not linked, by Lemma \ref{OHK4}
and Lemma \ref{evenK4}, $G[1,2,3,4]$ contains two $1$-homologous cycles. Without
loss
of generality, let $(1,2,4)$ and $(1,3,4)$ be $1$-homologous cycles. 

\vspace{\baselineskip}

Similarly, $G[2,4,5,6]$ contains two $1$-homologous cycles. The cycle $(4,5,6)$
is $0$-homologous by assumption and since $(2,5,6)$ is disjoint from $(1,3,4)$,
which is $1$-homologous, $(2,5,6)$ is $0$-homologous since $G$ is assumed to have no non-split link. Thus, $(2,4,5)$ and
$(2,4,6)$ are $1$-homologous
cycles. 

\vspace{\baselineskip}

In addition, $G[1,2,3,6]$ contains two $1$-homologous cycles. Since $(1,2,3)$
is $0$-homologous by assumption and $(1,3,6)$ is disjoint from $(2,4,5)$, which
is $1$-homologous, $(1,2,6)$ and $(2,3,6)$ are $1$-homologous. 

\vspace{\baselineskip}

Finally, $G[1,3,5,6]$ contains two $1$-homologous cycles. But, $(2,4,6)$,
$(2,4,5)$, and $(1,2,4)$ are $1$-homologous and disjoint from $(1,3,5)$,
$(1,3,6)$, and $(3,5,6)$, respectively, a contradiction, since $G[1,3,5,6]$
must contain two $1$-homologous cycles. 
\end{proof}

\begin{figure}
\begin{center}
\vskip -1.5in
\includegraphics[scale=0.5]{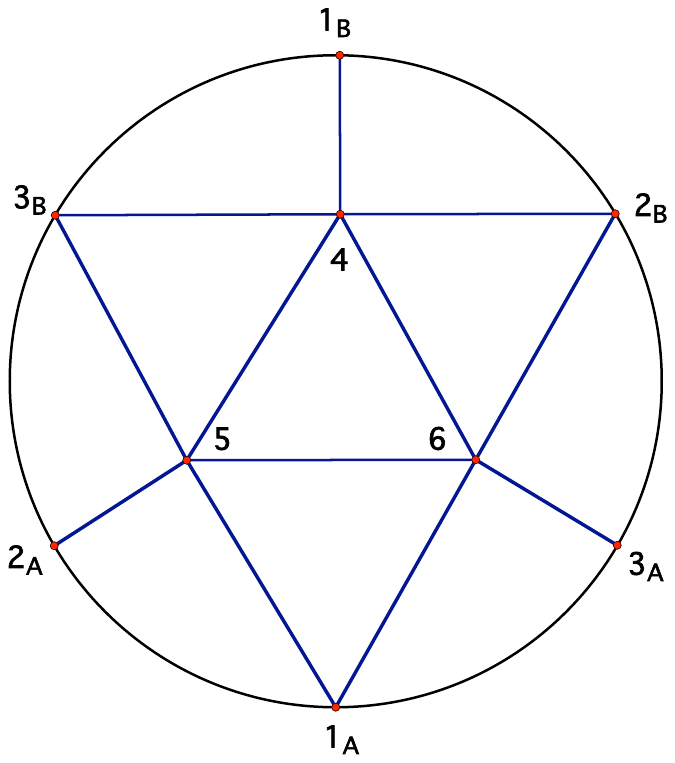}
\vskip -1.5in
\caption{A projection of a linkless embedding of $K_6$ in $\mathbb{R}P^3$.}
\label{linkless}
\end{center}
\end{figure}

\begin{prop} \label{K6}
 {\it Up to ambient isotopy and crossing changes, Figure
\ref{linkless} describes the only way to linklessly embed $K_6$ in
$\mathbb{R}P^3$.}
\end{prop}

\begin{proof}
Let $G$ be a complete graph on vertex set $\{1,2,3,4,5,6\}$. Embed $G$ in $\R
P^3$ linklessly. The graph $G$ contains a $0$-homologous $3$-cycle, since
otherwise $G$ contains two disjoint $1$-homologous cycles and is linked. Without
loss of generality, let $(4,5,6)$ be a $0$-homologous $3$-cycle. Proposition
\ref{TD} implies that the cycle $(1,2,3)$ is $1$-homologous as it is disjoint
from $(4,5,6)$. 

\vspace{\baselineskip}

Mroczkowski \cite{M} showed that every cycle can be made, via crossing changes
and ambient isotopy, into an unknotted $0$-cycle or the $1$-homologous cycle $g$ as
explained in the Introduction. Apply crossing changes and ambient
isotopy so that the embedding has a projection with vertices as drawn in Figure
\ref{linkless}. A priori, the edges between vertices $\{1,2,3\}$ and
$\{4,5,6\}$ may be more complicated than as drawn in the figure. 

\vspace{\baselineskip}

Vertices $\{1,2,3\}$ and the $3$-cycle $(1,2,3)$ lie on the boundary. In the
projection, we label the pair of antipodal identified vertices by $\{v_A,
v_B\}$ for $v \in \{1,2,3\}$.

\vspace{\baselineskip}

Consider the edge $E$ between $1$ 
and $4$.  Together with the path $(1_B,4)$ pictured in Figure
\ref{linkless}, 
it forms either a $0$-homologous or a $1$-homologous cycle.  If the 
cycle formed is $0$-homologous, then by Mroczowski's result, $E \cup (1_B, 4)$ can 
be made into the unknot by crossing changes, and then deformed so that 
$E$ is within a small neighborhood of the path $(1_B,4)$.  That is, the cycle does 
not cross the boundary of $D^2$.  If $E \cup (1_B,4)$ forms a $1$-homologous
cycle, then $E$ and the path formed by connecting $4$ to 
$1_A$ by a straight line segment form a $0$-homologous cycle.  By 
similar reasoning, the edge $E$ can be deformed, by crossing changes 
and ambient isotopy, to be within a small neighborhood of $(1_A,4)$; 
that is to say, it does not cross the boundary of $D^2$. By similar reasoning,
all edges between vertices $\{1,2,3\}$ and $\{4,5,6\}$ may be drawn in the
projection onto $\R P^2$ without crossing the boundary. 

\vspace{\baselineskip}

 We now describe how vertices $\{1,2,3\}$ connect to
vertices $\{4,5,6\}$. We use that $G$ does not contain two disjoint
$1$-homologous cycles or a $0$-homologous $K_4$ by Lemma \ref{OHK4}. 

\vspace{\baselineskip}

Let $v \in \{1,2,3\}$. Then $v$ connects to one of $\{4,5,6\}$ from $v_A$ and
connects to one of $\{4,5,6\}$ from $v_B$, otherwise, $G$ has a $0$-homologous
$K_4$ subgraph. Without loss of generality, suppose that $\{2_A\}$ connects to
$\{5\}$ and $\{2_B\}$ connects to $\{4\}$ and $\{6\}$. 

\vspace{\baselineskip}

If $\{1_A\}$ or $\{1_B\}$ connect to both $\{4\}$ and $\{6\}$, then
$G[1,2,4,6]$ is a $0$-homologous $K_4$. Thus, without loss of generality, let
$\{1_B\}$ connect to $\{4\}$ and $\{1_A\}$ connect to $\{6\}$. 

\vspace{\baselineskip}

Vertex $\{1_A\}$ connects to $\{5\}$ since otherwise, any arrangement of edges
connecting vertex $\{3\}$ to vertices $\{4,5,6\}$ induces either two disjoint
$1$-homologous cycles or a $0$-homologous $K_4$ subgraph, as shown in the table
below. 

\begin{center}
\begin{tabular}{c|c||c}
 Vertices $\{3_A\}$&Vertices $\{3_B\}$  & $1$-homologous cycles \\
 connects to& connects to& or $0$-homologous $K_4$ \\
\hline
$\{4\}$ & $\{5\},\{ 6\}$  &$(1,3,6)$, $(2,4,5)$\\
$\{5\}$ & $\{4\},\{ 6\}$ & $G[2,3,4,6]$ \\
$\{6\}$ & $\{4\},\{5\}$ & $G[1,3,4,5]$ \\
$\{4\},\{ 5\}$ & $\{6\}$ & $(1,3,6)$, $(2,4,5)$ \\
$\{4\}, \{ 6\}$ & $\{5\}$ & $G[2,3,4,6]$ \\
$\{5\}, \{ 6\}$ & $\{4\}$ & $(1,2,5)$, $(3,4,6)$ 
\end{tabular}
\end{center}

Finally, the following table shows that vertex $\{3_A\}$ connects to $\{6\}$ and
vertex $\{3_B\}$ connects to $\{4,5\}$. Indeed, all other arrangements lead to
either two disjoint $1$-homologous cycles or a $0$-homologous $K_4$ subgraph. 

\begin{center}
\begin{tabular}{c|c||c}
Vertices $\{3_A\}$  & Vertices $\{3_B\}$  & $1$-homologous cycles \\
connects to& connects to& or $0$-homologous $K_4$ \\
\hline
$\{4\},\{ 5\}$ & $\{6\}$ & $(1,3,6)$, $(2,4,5)$ \\
$\{4\}, \{ 6\}$ & $\{5\}$ & $G[2,3,4,6]$ \\
$\{5\}, \{ 6\}$ & $\{4\}$ & $G[1,3,5,6]$ \\
$\{4\}$ & $\{5\},\{ 6\}$  &$G[1,3,4,6]$\\
$\{5\}$ & $\{4\},\{ 6\}$ & $G[2,3,4,6]$ \\
\end{tabular}
\end{center}
Thus, up to crossing changes and ambient isotopy, Figure \ref{linkless}
 depicts the
only
way $K_6$ may be linklessly embedded in $\R P^3$. 
\end{proof}

\begin{figure}
\begin{center}
\vskip -1.5in
\includegraphics[scale=0.45]{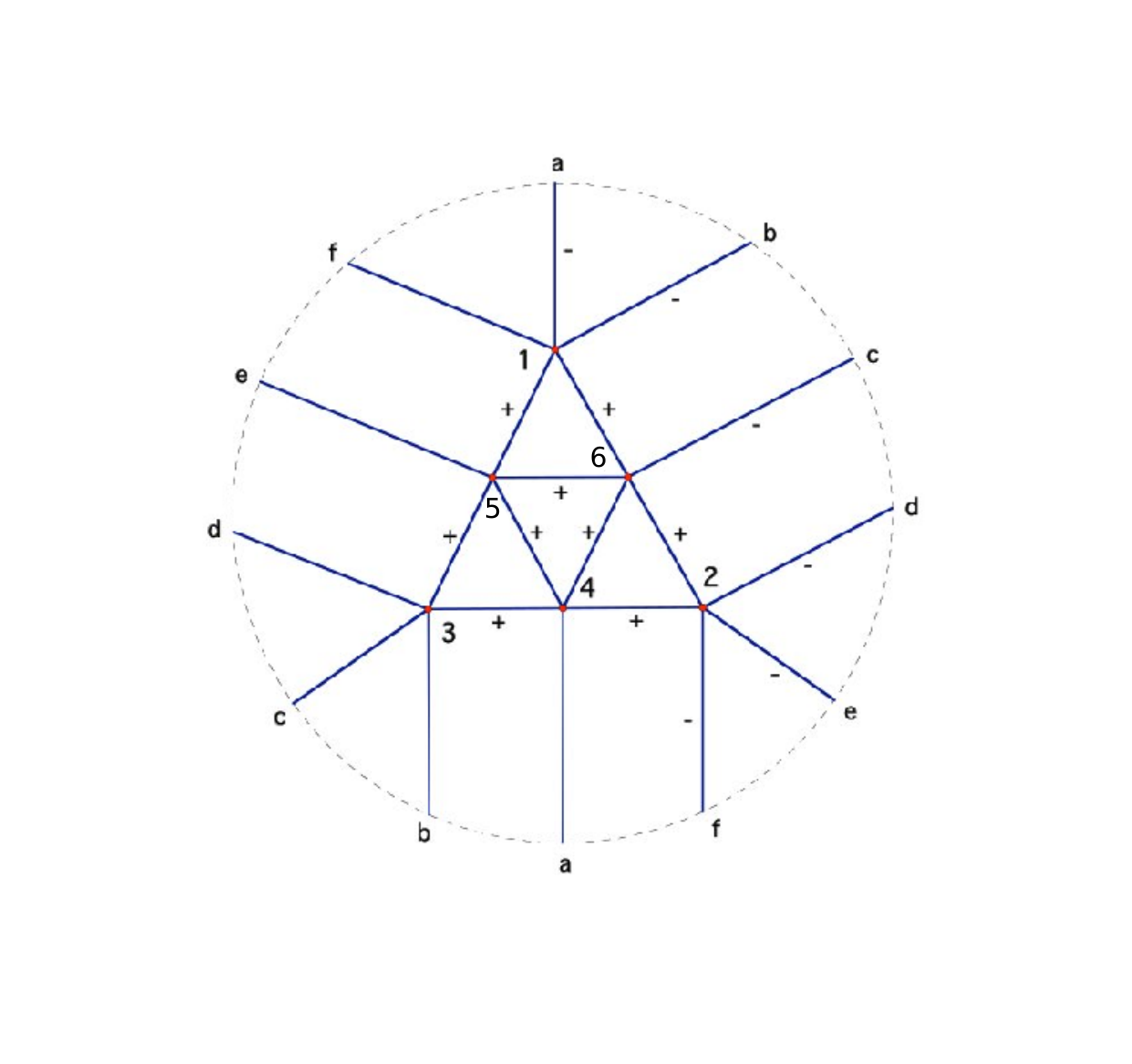}
\vskip -.5in
\caption{A signed linkless embedding of $K_6$ in $\mathbb{R}P^3$.}
\label{signed}
\end{center}
\end{figure}

Introduced by Harary in \cite{H}, {\it signed graphs} are graphs with each edge
assigned a $+$ or a $-$ sign, and constitute the final tool in our proof that
$K_{10}$ is intrinsically triple-linked in $\R P^3$. An embedding of a graph
$G$ into $\R P^3$ induces a signed graph as follows: deform the embedding to
that no vertices touch the bounding sphere in the model of $\R P^3$
with $\partial(D^3) \cong S^2$ and so that all intersections of edges
with the bounding sphere are transverse. Assign {\it $+$ edges} to be edges that
intersect the boundary an even number of times and {\it $-$ edges} to be edges
that intersect the boundary an odd number of times. An example is drawn in Figure \ref{signed}. Note that a cycle with an
odd number
of $-$ edges is $1$-homologous. 

\vspace{\baselineskip}

Two embeddings $G_1$ and $G_2$ of a graph $G$ are {\it crossing-change
equivalent} if $G_1$ can be obtained from $G_2$ by crossing changes and ambient
isotopy. By Proposition \ref{K6}, a linkless $K_6$ embedded in $\R P^3$ is
crossing-change equivalent to the embedding drawn in Figure \ref{signed}. That
is, if $G$ is a signed graph isomorphic to $K_6$ with vertex set
$\{1,2,3,4,5,6\}$, then $G$ is crossing-change equivalent to a signed graph
with $-$ edge set $S = \{(1,2), (1,3), (2,3), (1,4), (2,5), (3,6)\}$ and $+$
edge set $E(G)\backslash S$. 

\vspace{\baselineskip}

Our next result shows that if $G$ is a graph isomorphic to $K_{10}$, then $G$
is intrinsically triple-linked in $\R P^3$. We first sketch an outline. Using
results of \cite{FNP} and \cite{BF}, we show a $3$-linkless embedding of $G$,
if such an embedding exists, must contain a linkless $K_6$ subgraph. We prove the remaining four vertices must induce a $0$-homologous $K_4$ subgraph or the embedded graph contains a non-split triple-link. Finally, we determine the signs of the edges connecting the $K_6$ subgraph to the $K_4$ subgraph, eventually determining that any possible sign assignment results in a triple-link. Thus, no $3$-linkless embedding of $G$ can exist. 

\begin{thm} {\it The graph $K_{10}$ is intrinsically triple-linked in
$\mathbb{R}P^3$.}
\end{thm}

\begin{proof}

Let $G$ be a graph isomorphic to $K_{10}$ with vertex set
$\{1,2,3,4,5,6,7,8,9,10\}$. Embed $G$ in $\R P^3$ as a signed graph and assume,
toward a
contradiction, that $G$ is $3$-linkless. 

\vspace{\baselineskip}

If every subgraph of $G$ isomorphic to $K_6$ is linked, then Flapan, Naimi, and
Pommersheim's proof \cite{FNP} that $K_{10}$ is intrinsically linked in $\R^3$
nearly holds in $\R P^3$. However, at the end of their proof, they use that
$K_{3,3,1}$ is intrinsically linked in $\R^3$, but this graph embeds linklessly
in $\R P^3$. Bowlin and Foisy, \cite{BF}, modify the proof in \cite{FNP} to
only use the fact that $K_6$ is intrinsically linked in $\R^3$. Thus, in the
case that every subgraph of $G$ isomorphic to $K_6$ is linked, $G$ contains a
triple-link. So, we may assume that there exists a linkless $K_6$ subgraph of
$G$. Without loss of generality, suppose that $G[1,2,3,4,5,6]$ is linkless. By
Proposition \ref{K6}, $G[1,2,3,4,5,6]$ has an embedding that is crossing-change
equivalent to that drawn in Figure \ref{signed}. In particular, since crossing changes do not change the homology of cycles, we may assume $(1,2,3)$ is $1$-homologous. 

\vspace{\baselineskip}

\noindent {\it Claim:} The embedded induced subgraph $G[7,8,9,10]$ is
$0$-homologous. 

\begin{proof}
 Suppose $G[7,8,9,10]$ has a $1$-homologous cycle. Without loss of
generality, suppose $(7,8,9)$ is $1$-homologous. If
$G[4, 5, 6, 10]$ is not $0$-homologous, then two of $(4, 5, 10)$, $(4, 6, 10)$,
and $(5, 6, 10)$ are $1$-homologous by Lemma \ref{evenK4}, since we have
assumed $(4,5, 6)$ is $0$-homologous. Then $(1, 2, 3)$, $(7, 8, 9)$, and a cycle
from $G[4, 5, 6, 10]$ comprise three disjoint $1$-homologous cycles, so $G$ is
triple-linked. Thus, $G[4, 5, 6, 10]$ is $0$-homologous and so $G[1, 2, 4, 5, 6,
10]$ has a pair of linked cycles by Lemma \ref{OHK4}. Since $(7, 8, 9)$ is
$1$-homologous, and $(7, 8, 9)$ is disjoint from all the $1$-homologous cycles
in the second column of Table 1, Lemma \ref{LL} applies and $G$ has a
triple-link. Thus, $G[7, 8, 9, 10]$ is $0$-homologous. \end{proof}

\vspace{\baselineskip}

\begin{center}
\begin{tabular}{|c|c|}
\hline
Possible Linked& $1$-Homologous Cycle that  \\
Cycles in $G[1, 2, 4, 5, 6, 10]$ & shares an edge with a linked cycle \\
\hline
$(1, 2, 4)$, $(5, 6, 10)$ & $(1, 2, 3)$\\
$(1, 2, 5)$, $(4, 6, 10)$ & $(1, 2, 3)$\\
$(1, 2, 6)$, $(4, 5, 10)$ & $(1, 2, 3)$\\
$(1, 2, 10)$, $(4, 5, 6)$ & $(1, 2, 3)$\\
$(1, 4, 5)$, $(2, 6, 10)$ & $(1, 3, 5)$\\
$(1, 4, 6)$, $(2, 5, 10)$ & $(1, 4, 6)$\\
$(1, 4, 10)$, $(2, 5, 6)$ & $(2, 5, 6)$\\
$(1, 5, 6)$, $(2, 4, 10)$ & $(1, 3, 5)$\\
$(1, 5, 10)$, $(2, 4, 6)$ & $(1, 3, 5)$\\
$(1, 6, 10)$, $(2, 4, 5)$ & $(2, 4, 5)$\\
\hline
\end{tabular}

\vspace{.5\baselineskip}
Table $1$.
\end{center}

Since ambient isotopy and crossing changes do not change the homology of
cycles, we may modify the embedding of $G$ so that all edges in $G[7,8,9,10]$
are $+$ edges and the edges in $G[1,2,3,4,5,6]$ are $+$ and $-$ edges as
defined in Figure \ref{signed}. Many of the remaining arguments rely on
linked $K_6$ subgraphs of $G$ and use the argument highlighted in Table $1$. In
particular, though the $K_6$ subgraph of the modified embedding may contain a
different pair of linked cycles than the original embedding, our argument
relies only on the existence of linked cycles, not on the specific pair of
linked cycles. Thus, we now consider the signs of the edges connecting
$G[1,2,3,4,5,6]$ to $G[7,8,9,10]$. 

\vspace{\baselineskip}

\noindent {\it Claim:} If $v \in \{1,2,3\}$, then edges from $v$ to
$G[7,8,9,10]$ have the same sign. 

\begin{proof}
Assume toward a contradiction that the edges from $\{1\}$ to $G[7,8,9,10]$ do
not all have the same sign. Without loss of generality, let $(1,7)$ be a $+$
edge and $(1,8)$ a $-$ edge. Then $(1,7,8)$ is a $1$-homologous cycle. 

\vspace{\baselineskip}

Consider $G[3,4,6,9]$. Since $(3,4,6)$ is $1$-homologous, $G[3,4,6,9]$ contains
another $1$-homologous cycle by Lemma \ref{evenK4}. If $(3,4,9)$ or
$(3,6,9)$ is $1$-homologous then the sets $\{(1,7,8), (2,5,6), (3,4,9)\}$ or
$\{(1,7,8), (2,4,5), (3,6,9)\}$ form three disjoint $1$-homologous cycles,
respectively, and so $G$ is triple-linked. Thus, $(4,6,9)$ is the second
$1$-homologous cycle in $G[3,4,6,9]$. 

\vspace{\baselineskip}

Since $(2,3,4)$ is $1$-homologous, the induced subgraph $G[2,3,4,9]$ contains a
second $1$-homologous cycle by Lemma \ref{evenK4}. As shown above, $(3,4,9)$ is
$0$-homologous. If $(2,4,9)$ is $1$-homologous, then $(1,7,8)$, $(2,4,9)$, and
$(3,5,6)$ form three disjoint $1$-homologous cycles, so $G$ is triple-linked.
So, $(2,3,9)$ is the second $1$-homologous cycle in $G[2,3,4,9]$. 

\vspace{\baselineskip}

Similarly, since $(3,5,6)$ is $1$-homologous, $G[3,5,6,9]$ contains a second
$1$-homologous cycle by Lemma \ref{evenK4}. As shown above, $(3,6,9)$ is
$0$-homologous. Additionally, $(5,6,9)$ is $0$-homologous, otherwise $(1,7,8),
(2,3,4)$, and $(5,6,9)$ form three disjoint $1$-homologous cycles and $G$ is
triple-linked. Thus, $(3,5,9)$ is a $1$-homologous cycle. 

\vspace{\baselineskip}

As $(1,7,8)$ and $(4,6,9)$ are $1$-homologous, $G[2,3,5,10]$ is a
$0$-homologous $K_4$, since, otherwise, $G$ contains three disjoint
$1$-homologous cycles. By Lemma \ref{OHK4}, $G[2,3,4,5,6,10]$ contains a pair
of linked cycles. Since $(1,7,8)$ is $1$-homologous and disjoint from all of
the $1$-homologous cycles in the second column of Table $2$, Lemma \ref{LL}
applies and $G$ contains a triple-link, a contradiction. 

\vspace{\baselineskip}

\begin{center}
\begin{tabular}{|c|c|}
\hline
Possible Linked& $1$-Homologous Cycle that  \\
Cycles in $G[2, 3, 4, 5, 6, 10]$ & shares an edge with a linked cycle \\
\hline
$(2, 3, 4)$, $(5, 6, 10)$ & $(2, 3, 4)$\\
$(2, 3, 5)$, $(4, 6, 10)$ & $(4, 6, 9)$\\
$(2, 3, 6)$, $(4, 5, 10)$ & $(2, 3, 9)$\\
$(2, 3, 10)$, $(4, 5, 6)$ & $(2, 3, 9)$\\
$(2, 4, 5)$, $(3, 6, 10)$ & $(2, 4, 5)$\\
$(2, 4, 6)$, $(3, 5, 10)$ & $(4, 6, 9)$\\
$(2, 4, 10)$, $(3, 5, 6)$ & $(3, 5, 9)$\\
$(2, 5, 6)$, $(3, 4, 10)$ & $(2, 5, 6)$\\
$(2, 5, 10)$, $(3, 4, 6)$ & $(4, 6, 9)$\\
$(2, 6, 10)$, $(3, 4, 5)$ & $(3, 5, 9)$\\
\hline
\end{tabular}

\vspace{.5\baselineskip}
Table $2$.
\end{center}

\vspace{\baselineskip}

Thus, $\{1\}$ connects to $G[7,8,9,10]$ via all $+$ edges or all $-$ edges, and
similar reasoning applies to vertices $\{2\}$ and $\{3\}$. \end{proof}

A similar argument, using different induced subgraphs, show the edges between
each of the remaining vertices of $G[1,2,3,4,5,6]$ and $G[7,8,9,10]$ also have
the same sign. 

\vspace{\baselineskip}

\noindent {\it Claim:} If $v \in \{4,5,6\}$, then all edges from $v$ to
$G[7,8,9,10]$ have the same sign. 

\begin{proof}
Towards a contradiction, suppose that not all the edges from $\{4\}$ to
$G[7,8,9,10]$ have the same sign. Without loss of generality, let $(4,7)$ be a
$+$ edge and $(4,8)$ be a $-$ edge. Then $(4,7,8)$ is a $1$-homologous cycle. 

\vspace{\baselineskip}

Since $(1,2,3)$ is a $1$-homologous cycle, $G[1,2,3,9]$ contains a second
$1$-homologous cycle by Lemma \ref{evenK4}. If $(1,3,9)$ or $(1,2,9)$ are
$1$-homologous, then $\{(1,3,9), (2,5,6), (4,7,8)\}$ or $\{(1,2,9), (3,5,6),
(4,7,8)\}$ form three disjoint $1$-homologous cycles, respectively. So,
$(2,3,9)$ is the second $1$-homologous cycle in $G[1,2,3,9]$. 

\vspace{\baselineskip}

Since $(2,3,9)$ and $(4,7,8)$ are $1$-homologous, $G[1,5,6,10]$ is a
$0$-homologous $K_4$, otherwise, $G$ contains three disjoint $1$-homologous
cycles. By Lemma \ref{OHK4}, $G[1,2,3,5,6,10]$ contains a pair of linked
cycles. Since $(4,7,8)$ is $1$-homologous and disjoint from all $1$-homologous
cycles in the second column of Table $3$, Lemma \ref{LL} applies and $G$
contains a triple link. 

\vspace{\baselineskip}

\begin{center}
\begin{tabular}{|c|c|}
\hline
Possible Linked& $1$-Homologous Cycle that  \\
Cycles in $G[1, 2, 3, 5, 6, 10]$ & shares an edge with a linked cycle \\
\hline
$(1, 2, 3)$, $(5, 6, 10)$ & $(1, 2, 3)$\\
$(1, 2, 5)$, $(3, 6, 10)$ & $(2, 5, 9)$\\
$(1, 2, 6)$, $(3, 5, 10)$ & $(1, 2, 6)$\\
$(1, 2, 10)$, $(3, 5, 6)$ & $(3, 5, 6)$\\
$(1, 3, 5)$, $(2, 6, 10)$ & $(1, 3, 5)$\\
$(1, 3, 6)$, $(2, 5, 10)$ & $(2, 5, 9)$\\
$(1, 3, 10)$, $(2, 5, 6)$ & $(2, 5, 9)$\\
$(1, 5, 6)$, $(2, 3, 10)$ & $(2, 3, 9)$\\
$(1, 5, 10)$, $(2, 3, 6)$ & $(2, 3, 9)$\\
$(1, 6, 10)$, $(2, 3, 5)$ & $(2, 3, 9)$\\
\hline
\end{tabular}

\vspace{.5\baselineskip}
Table $3$.
\end{center}

\vspace{\baselineskip}

Therefore, all edges from $\{4\}$ to $G[7,8,9,10]$ have the same sign. A
similar argument show that all edges from vertices $\{5\}$ and $\{6\}$ to
$G[7,8,9,10]$ have the same sign. \end{proof}

The previous two claims show that the edges from each vertex in
$G[1,2,3,4,5,6]$ to the vertices of $G[7,8,9,10]$ have the same sign. As we
have assigned signs to the edges of $G[1,2,3,4,5,6]$ and $G[7,8,9,10]$, there
remain $2^6$ possible embedding classes. We consider all cases. If all edges
from vertex $v \in \{1,2,3,4,5,6\}$ to $G[7,8,9,10]$ are $+$ edges, we write
$v_{+}$, and otherwise $v_{-}$. For $v_{x}$ with $x \in \{+, -\}$, we say ``the
sign of vertex $v$ is $x$.''

\vspace{\baselineskip}

\noindent {\it Claim:} The two vertices in each of the pairs $\{1,4\}$,
$\{2,5\}$, and $\{3,6\}$ have different signs. 

\begin{proof}
Suppose toward a contradiction that $\{1\}$ and $\{4\}$ are both $+$ edges.
Then $(1,4,7)$ is a $1$-homologous cycle. 

\vspace{\baselineskip}

Since both $(2,5)$ and $(3,6)$ are $-$ edges, if both pairs of vertices
$\{2,5\}$ and $\{3,6\}$ share the same sign (eg. $2_{+}, 5_{+}, 3_{-}, 6_{-}$),
then $(2,5,8)$ and $(3,6,9)$ are $1$-homologous cycles. Then, $(1,4,7)$,
$(2,5,8)$, and $(3,6,9)$ are disjoint $1$-homologous cycles, so $G$ is
triple-linked. 

\vspace{\baselineskip}

Since both $(2,6)$ and $(3,5)$ are $+$ edges, if both pairs of vertices
$\{2,6\}$ and $\{3,5\}$ have different signs (eg. $2_{+}, 6_{-}, 3_{+}, 5_{-}$),
then $(2,6,8)$ and $(3,5,9)$ are $1$-homologous cycles. Then, $(1,4,7)$,
$(2,6,8)$, and $(3,5,9)$ are disjoint $1$-homologous cycles, so $G$ is
triple-linked. 

\vspace{\baselineskip}

The edge $(2,3)$ is a $-$ edge and $(5,6)$ is a $+$ edge, so if $\{2\}$ and
$\{3\}$ share the same sign and $\{5\}$ and $\{6\}$ have different signs, (eg.
$2_{+}, 3_{+}, 5_{+}, 6_{-}$), then $(2,3,8)$ and $(5,6,9)$ are $1$-homologous
cycles. So, $(1,4,7), (2,3,8)$, and $(5,6,9)$ form disjoint $1$-homologous
cycles, so $G$ is triple-linked. 

\vspace{\baselineskip}

If $G$ is embedded with either $\{2_{-}, 3_{+}, 5_{+}, 6_{-} \}$ or $\{ 2_{+},
3_{-}, 5_{-}, 6_{+}\}$, then $(1,4,7)$ and $(5,6,8)$ are disjoint
$1$-homologous cycles, so $G[2,3,9,10]$ is a $0$-homologous $K_4$, or $G$ has
a triple-link. So, by Lemma \ref{OHK4}, $G[1,2,3,4,9,10]$ has a pair of linked
cycles. Since $(5,6,8)$ is $1$-homologous and is disjoint from all of
the $1$-homologous cycles in the second column of Table $4$, $G$ has a triple
link by Lemma \ref{LL}. 

\vspace{\baselineskip}

\begin{center}
\begin{tabular}{|c|c|}
\hline
Possible Linked& $1$-Homologous Cycle that  \\
Cycles in $G[1, 2, 3, 4, 9, 10]$ & shares an edge with a linked cycle \\
\hline
$(1, 2, 3)$, $(4, 9, 10)$ & $(1, 2, 3)$\\
$(1, 2, 4)$, $(3, 9, 10)$ & $(1, 4, 7)$\\
$(1, 2, 9)$, $(3, 4, 10)$ & $(1, 2, 7)$\\
$(1, 2, 10)$, $(3, 4, 9)$ & $(1, 2, 7)$\\
$(1, 3, 4)$, $(2, 9, 10)$ & $(1, 4, 7)$\\
$(1, 3, 9)$, $(2, 4, 10)$ & $(1, 3, 7)$\\
$(1, 3, 10)$, $(2, 4, 9)$ & $(1, 3, 7)$\\
$(1, 4, 9)$, $(2, 3, 10)$ & $(1, 4, 7)$\\
$(1, 4, 10)$, $(2, 3, 9)$ & $(1, 4, 7)$\\
$(1, 9, 10)$, $(2, 3, 4)$ & $(2, 3, 4)$\\
\hline
\end{tabular}

Table $4$.
\end{center}

Finally, if $G$ is embedded with one of the remaining configurations,
$\{\{2_{-}, 3_{+}, 5_{+}, 6_{+} \}$, $\{2_{-}, 3_{+}, 5_{-}, 6_{-} \},\{2_{+},
3_{-}, 5_{-}, 6_{-} \},  \{2_{-}, 3_{+}, 5_{+}, 6_{+} \}   \}$, then one of
$\{(2,5,6), (3,5,6)\}$ is $1$-homologous. Since $G[7,8,9,10]$ is a
$0$-homologous $K_4$, $G[1,4,7,8,9,10]$ contains a pair of linked cycles by
Lemma \ref{OHK4}. Both $(2,5,6)$ and $(3,5,6)$ are disjoint from one
$1$-homologous cycle in each row of the second column of Table $5$. Thus, by
Lemma \ref{LL}, $G$ is triple-linked. 

\vspace{\baselineskip}

\begin{center}
\begin{tabular}{|c|c|}
\hline
Possible Linked& $1$-Homologous Cycle that  \\
Cycles in $G[1, 4, 7, 8, 9, 10]$ & shares an edge with a linked cycle \\
\hline
$(1, 4, 7)$, $(8, 9, 10)$ & $(1, 4, 7)$\\
$(1, 4, 8)$, $(7, 9, 10)$ & $(1, 4, 8)$\\
$(1, 4, 9)$, $(7, 8, 10)$ & $(1, 4, 9)$\\
$(1, 4, 10)$, $(7, 8, 9)$ & $(1, 4, 10)$\\
$(1, 7, 8)$, $(4, 9, 10)$ & $(1, 2, 7), (1,3,7)$\\
$(1, 7, 9)$, $(4, 8, 10)$ & $(1, 2, 7), (1,3,7) $\\
$(1, 7, 10)$, $(4, 8, 9)$ & $(1, 2, 7), (1,3,7)$\\
$(1, 8, 9)$, $(4, 7, 10)$ & $(1, 2, 8), (1,3,8)$\\
$(1, 8, 10)$, $(4, 7, 9)$ & $(1, 2, 8), (1,3,8)$\\
$(1, 9, 10)$, $(4, 7, 8)$ & $(1, 2, 9), (1,3,9)$\\
\hline
\end{tabular}

\vspace{.5\baselineskip}
Table $5$.
\end{center}

\vspace{\baselineskip}

So, in each embedding of $G$ with $1_{+}$ and $4_{+}$, $G$ contains a triple
link. A similar argument holds in the case that $G$ is embedded with $1_{-}$
and $4_{-}$ and for the other vertex pairs $\{2,5\}$ and $\{3,6\}$. 
\end{proof}

We now suppose $G$ is embedded with $1_{+}$ and $4_{-}$. By the last claim, the
vertices in each of the pairs $\{2,5\}$ and $\{3,6\}$ have different signs. So,
there are four cases to consider: $\{ \{2_{+}, 3_{+}, 5_{-}, 6_{-}\}, \{2_{+},
3_{-}, 5_{-}, 6_{+}\},\{2_{-}, 3_{+}, 5_{+}, 6_{-}\},\{2_{-}, 3_{-}, 5_{+},
6_{+}\}\}$. 

\vspace{\baselineskip}

First, if the embedding has $\{2_{+}, 3_{+}, 5_{-}, 6_{-}\}$, then $(1,6,7)$,
$(2,4,9)$, and $(3,5,8)$ form three disjoint $1$-homologous cycles, so $G$ is
triple-linked. Second, suppose the embedding has $\{2_{+},
3_{-}, 5_{-}, 6_{+}\}$ or $\{2_{-}, 3_{+}, 5_{+}, 6_{-}\}$. Then the second
column of Table $6$ contains $1$-homologous cycles. Since $G[7,8,9,10]$ is
$0$-homologous, $G[4,6,7,8,9,10]$ has a pair of linked cycles by Lemma
\ref{OHK4}. Since $(1,2,3)$ is $1$-homologous and disjoint from all
$1$-homologous cycles in the second column of Table $6$, Lemma \ref{LL} applies
and $G$ contains a triple link. 

\vspace{\baselineskip}

\begin{center}
\begin{tabular}{|c|c|}
\hline
Possible Linked& $1$-Homologous Cycle that  \\
Cycles in $G[4, 6, 7, 8, 9, 10]$ & shares an edge with a linked cycle \\
\hline
$(4, 6, 7)$, $(8, 9, 10)$ & $(4, 6, 7)$\\
$(4, 6, 8)$, $(7, 9, 10)$ & $(4, 6, 8)$\\
$(4, 6, 9)$, $(7, 8, 10)$ & $(4, 6, 9)$\\
$(4, 6, 10)$, $(7, 8, 9)$ & $(4, 6, 10)$\\
$(4, 7, 8)$, $(6, 9, 10)$ & $(5, 6, 9)$\\
$(4, 7, 9)$, $(6, 8, 10)$ & $(5, 6, 8)$\\
$(4, 7, 10)$, $(6, 8, 9)$ & $(5, 6, 8)$\\
$(4, 8, 9)$, $(6, 7, 10)$ & $(5, 6, 7)$\\
$(4, 8, 10)$, $(6, 7, 9)$ & $(5, 6, 7)$\\
$(4, 9, 10)$, $(6, 7, 8)$ & $(5, 6, 7)$\\
\hline
\end{tabular}

\vspace{.5\baselineskip}
Table $6$.
\end{center}

\vspace{\baselineskip}

Finally, if the embedding has $\{2_{-}, 3_{-}, 5_{+},6_{+}\}$, then the second
column of Table $7$ contains $1$-homologous cycles. As above, since $(1,2,3)$ is
$1$-homologous and disjoint from all $1$-homologous cycles in the second column
of Table $7$, Lemma \ref{LL} applies and $G$ contains a triple link.

\vspace{\baselineskip}

\begin{center}
\begin{tabular}{|c|c|}
\hline
Possible Linked& $1$-Homologous Cycle that  \\
Cycles in $G[4, 6, 7, 8, 9, 10]$ & shares an edge with a linked cycle \\
\hline
$(4, 6, 7)$, $(8, 9, 10)$ & $(4, 6, 7)$\\
$(4, 6, 8)$, $(7, 9, 10)$ & $(4, 6, 8)$\\
$(4, 6, 9)$, $(7, 8, 10)$ & $(4, 6, 9)$\\
$(4, 6, 10)$, $(7, 8, 9)$ & $(4, 6, 10)$\\
$(4, 7, 8)$, $(6, 9, 10)$ & $(4, 5, 7)$\\
$(4, 7, 9)$, $(6, 8, 10)$ & $(4, 5, 7)$\\
$(4, 7, 10)$, $(6, 8, 9)$ & $(4, 5, 7)$\\
$(4, 8, 9)$, $(6, 7, 10)$ & $(4, 5, 8)$\\
$(4, 8, 10)$, $(6, 7, 9)$ & $(4, 5, 8)$\\
$(4, 9, 10)$, $(6, 7, 8)$ & $(4, 5, 9)$\\
\hline
\end{tabular}

\vspace{.5\baselineskip}
Table $7$.
\end{center}

\vspace{\baselineskip}

The same argument holds if $G$ is embedded with $1_{-}$ and $4_{+}$. So, for
any assignment of signs to the edges from $\{1\}$ and $\{4\}$ to
$G[7,8,9,10]$, $G$ contains a triple-link, a contradiction. Thus, every
embedding of $G$ into $\R P^3$ contains a triple-link, so $G$ is intrinsically
triple-linked in $\R P^3$. \end{proof}

Flapan, Naimi, and Pommersheim \cite{FNP} show that $K_9$ can be embedded
$3$-linklessly in $\R^3$, and so $K_9$ can be embedded $3$-linklessly in $\R
P^3$ as well. Thus, $10$ is the smallest $n$ for which $K_n$ is intrinsically
triple
linked in $\R P^3$.

\section{Other intrinsically triple-linked graphs in $\R P^3$}

In this section, we exhibit other intrinsically triple-linked graphs in
$\R P^3$. We show that two graphs shown in \cite{BF} to be intrinsically
triple-linked in $\R^3$ may be embedded $3$-linklessly in $\R P^3$. Moreover,
the graphs obtained by taking two disjoint copies of these graphs described in
\cite{BF} give intrinsically triple-linked graphs in $\R P^3$. We begin by
describing a family of intrinsically $n$-linked graphs in $\R P^3$. 

\begin{lemma}\label{cc}
 If an embedded graph has all $0$-homologous cycles, then it is crossing-change
equivalent to a spatial embedding. 
\end{lemma}
\begin{proof}
 Take a spanning tree in the embedded graph.  Since a 
spanning tree is contractible, it can be deformed so that none of its 
edges touch the boundary of $D^2$.  Order the edges that do not lie in 
the spanning tree.  Now take the first edge not in the spanning tree.  
If this edge does not touch the boundary, move on to the next edge.  
Otherwise, the edge lies in a cycle that, by assumption, is $0$-
homologous.  By Mroczowski's result, the cycle can be made into an 
unknot by crossing changes.  Since the unknot is $0$-homologous, it 
bounds a disk.  Deform the edge by pulling in the disk towards the 
edges of the cycle that lie in the spanning tree.  Thus, the edge can 
be deformed so that it does not touch the boundary of $D^2$.  
Eventually, all of the edges not in the spanning tree can be deformed, 
if necessary, not to touch the boundary.  The resulting embedding is 
equivalent to a spatial embedding.  Thus, the original embedding was 
crossing-change equivalent to a spatial embedding.
\end{proof}

\begin{prop} {\it A graph composed of $n$ disjoint copies of an intrinsically
$n$-linked graph in $\mathbb{R}^3$ is intrinsically $n$-linked in
$\mathbb{R}P^3$. In particular, three disjoint copies of intrinsically
triple-linked graphs in $\mathbb{R}^3$ are intrinsically triple-linked in
$\mathbb{R}P^3$.} \end{prop}

\begin{proof}
 Let $G$ be a graph that is intrinsically $n$-linked in $\R^3$, and let $G_i$
be isomorphic to $G$ for $i = 1, \ldots, n$. Let $\Gamma = \sqcup_{i=1}^n G_i$
be the disjoint union of $n$ graphs isomorphic to $G$. If $G_i$ contains all
$0$-homologous cycles for some $i$, then $G_i$ is crossing-change equivalent to
a spatial embedding by Lemma \ref{cc}. Thus, $G_i$, and hence $G$, is $n$-linked
in
$\R P^3$. 

\vspace{\baselineskip}

Otherwise, each $G_i$ contains a $1$-homologous cycle. Thus, $\Gamma$ contains
$n$ disjoint $1$-homologous cycles, so contains an $n$-link. Therefore,
$\Gamma$ is intrinsically $n$-linked in $\R P^3$. \end{proof}

The graph $K_{10}$ is an example of a one-component graph that is intrinsically
triple-linked in $\R P^3$. We now exhibit two intrinsically triple-linked graphs
in $\R P^3$, each comprised of two
components. In each case, the components are intrinsically
triple-linked in $\R^3$.
The question remains whether there exists a minor-minimal intrinsically
triple-linked graph of three components in $\R P^3$. 

\vspace{\baselineskip}

Bowlin and Foisy prove the following graphs are intrinsically linked in $\R^3$. 

\begin{thm} {\cite{BF} \label{PEPE} \it Let $G$ be a graph containing two
disjoint graphs
from the Petersen family, $G_1$ and $G_2$ as subgraphs.  If there are edges
between the two subgraphs $G_1$ and $G_2$ such that the edges form a 6-cycle
with vertices that alternate between $G_1$ and $G_2$, then $G$ is minor-minimal
intrinsically triple-linked in $\mathbb{R}^3$.}
\end{thm}

\begin{figure}
\begin{center}
\vskip -2.1in
\includegraphics[scale=0.6]{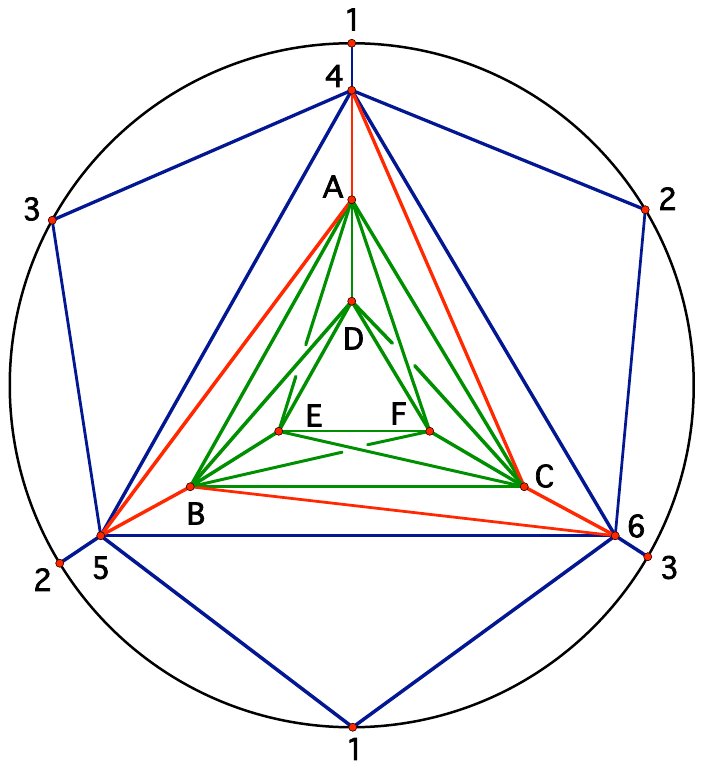}
\vskip -1.7in
\caption{A 3-linkless embedding of $K_6$ connected to $K_6$ along a 6-cycle in
$\mathbb{R}P^3$.}
\label{K6K6}
\end{center}
\end{figure}

If $G_1$ and $G_2$, as in the theorem, are isomorphic to $K_6$, this result
does not hold in $\R P^3$. A $3$-linkless embedding of $G = G_1 \sqcup G_2$ is
shown in Figure \ref{K6K6}. We now show that the graph obtained from two
disjoint
copies of $G$ is minor-minimal intrinsically triple-linked in $\R P^3$. 

\vspace{\baselineskip}

\begin{thm} \label{K6K6RP3}
 Let $G_1$ be a graph containing two disjoint copies of $K_6$ with edges
between the two $K_6$ subgraphs that form a $6$-cycle with vertices alternating
between the two $K_6$ subgraphs. If $G_2$ is a graph isomorphic to $G_1$ and $G
= G_1 \sqcup G_2$, then $G$ is minor-minimal intrinsically triple-linked in $\R
P^3$. 
\end{thm}

\begin{proof}
 Let $G = G_1 \sqcup G_2$ be as in the theorem, and embed $G$ in $\R P^3$.

\vspace{\baselineskip}

If either $G_1$ or $G_2$ contain all $0$-homologous cycles, then that subgraph
is crossing-change equivalent to a spatial embedding by Lemma \ref{cc}, and
hence
triple linked by Theorem \ref{PEPE}. Thus, $G$ contains a triple-link. So, now
suppose that both $G_1$ and $G_2$ contain a $1$-homologous cycle. 

\vspace{\baselineskip}

In both $G_1$ and $G_2$, any cycle of length greater than three can be
subdivided by an edge $e$ into a ``$\theta$-graph'': two cycles of smaller
length, disjoint, except for edge $e$. That is, there exists an edge $e = (v_1,
v_i)$ in $G[v_1, \ldots, v_n]$ so that $c = (v_1, \ldots, v_n)$ may be divided
into $c_1 \cup c_2 = (v_1, \ldots, v_i) \cup (v_i, \ldots, v_n , v_1)$. If $c$
is $1$-homologous, then in any signed embedding of $G$, $c$ has an odd number of
$-$ edges. So, either $c_1$ or $c_2$ has an odd number of $-$ edges, and is thus
$1$-homologous. By iterating this procedure, we conclude that both $G_1$ and
$G_2$ contain a $1$-homologous $3$-cycle. 

\vspace{\baselineskip}

Label the vertices of $G_1$ $\{1,2,3,4,5,6,A,B,C,D,E,F\}$ so that
$G[1,2,3,4,5,6]\cong K_6$ and $G[A,B,C,D,E,F] \cong K_6$ are connected by edges
$(4,A)$, $(4,C)$, $(5,A)$, $(5,B)$, $(6,B)$, and $(6,C)$. Up to isomorphism,
there are five $3$-cycle equivalence classes in $G_1$. The set $S = \{(1,2,3),
(1,2,4), (1,4,5), (4,5,6), (4,5,A)\}$ contains one representative from each
$3$-cycle equivalence class. So, without loss of generality, we may suppose
that $S$ contains a $1$-homologous $3$-cycle. 

\vspace{\baselineskip}

If $G[B,C,E,F] \cong K_4$ contains a $1$-homologous cycle, then this cycle,
along with the $1$-homologous cycle in $S$ and the $1$-homologous cycle in
$G_2$ form three disjoint $1$-homologous cycles and so $G$ contains a triple
link. Now suppose that $G[B,C,E,F]$ is $0$-homologous, so that $G[A,B,C,D,E,F]$
contains a pair of linked cycles by Lemma \ref{OHK4}. 

\vspace{\baselineskip}

First suppose that the $1$-homologous cycle, $c_1 \in S$ is not $(4,5,A)$. By
the pigeon hole principle, two vertices in $\{A,B,C\}$ are in one of the
components, $c_2$, of the linked cycles in $G[A,B,C,D,E,F]$. Use the edges of
the $6$-cycle to join $c_2$ to $c_1$ along disjoint paths. By Lemma \ref{BFLL},
$G$ contains a triple link. 

\vspace{\baselineskip}

Now suppose that the $1$-homologous cycle in $S$ is $(4,5,A)$. If there is a one
homologous cycle in $G[1, 2, 3, 6]$ then this cycle will link with $(4, 5, A)$
and the $1$-homologous cycle in $G_2$, so $G$ contains a triple-link. Else,
$G[1, 2, 3, 4, 5, 6]$ has a pair of linked cycles by Lemma \ref{OHK4}. By the
pigeon-hole principle, at least two vertices in the set $\{4, 5, 6\}$ are in a
linked cycle, $c_3$, within $G[1,2,3,4,5,6]$. Similarly, at least two vertices
of $\{A, B, C\}$ are in a linked cycle, $c_4$, within $G[A,B,C,D,E,F]$. As a
result of the $6$-cycle connecting these two copies of $K_6$, there are two
disjoint edges between $c_3$ and $c_4$. By Lemma \ref{BFLL}, $G$ is
triple-linked. 

\vspace{\baselineskip}

To see $G$ is minor-minimal with respect to intrinsic triple-linking in $\mathbb
{R}P^3$, embed $G$ so that $G_1$ is embedded as in the drawing in Figure
\ref{K6K6} and $G_2$ is contained in a sphere that lies in the complement of
$G_1$.  Therefore, $G_1$ does not have any triple-links and no cycle in $G_1$ is
linked with a cycle in $G_2$. Without loss of generality, if we delete an edge,
contract an edge or delete any vertex on $G_2$, it will have an affine linkless
embedding. Thus, we can re-embed $G_2$ within the sphere in each case.
Therefore, $G$
is minor-minimal for intrinsic triple-linking. \end{proof}

\vspace{\baselineskip}

Bowlin and Foisy prove the following graph is intrinsically triple-linked in
$\R^3$. 

\begin{thm} \cite{BF} \label{K7K7g} {\it Let $G$ be a graph formed by
identifying
an edge of
$K_7$ with an edge from another copy of $K_7$.  Then $G$ is intrinsically
triple-linked in $\mathbb{R}^3$.}
\end{thm}

\begin{figure}
\begin{center}
\vskip -1.5in
\includegraphics[scale=0.45]{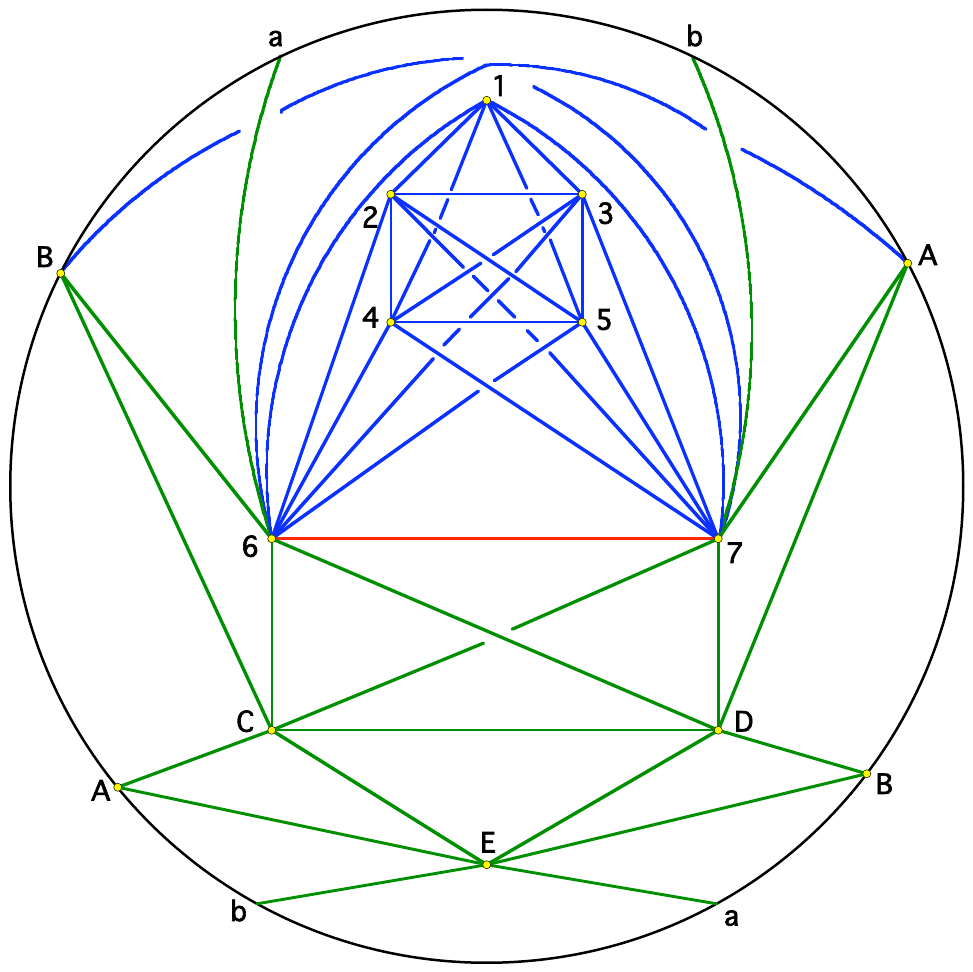}
\vskip -1in
\caption{A 3-linkless embedding of $K_7$ connected to $K_7$ along an edge in
$\mathbb{R}P^3$.}
\label{K7K7}
\end{center}
\end{figure}

The graph $G$ defined in Theorem \ref{K7K7g} may be embedded $3$-linklessly in
$\R P^3$, as drawn in Figure \ref{K7K7}. As in the previous result, the graph
consisting of two disjoint copies of this graph is intrinsically
linked in $\R P^3$. 

\begin{thm} \label{K7e}
 Let $G_1$ be a graph formed by identifying an edge of $K_7$ with an edge from
another copy of $K_7$. If $G_2$ is isomorphic to $G_1$ and $G= G_1 \sqcup G_2$
is the disjoint union of $G_1$ and $G_2$, then $G$ is intrinsically linked in
$\R P^3$. 
\end{thm}

\begin{proof}
Let $G = G_1 \sqcup G_2$ be as in the proposition, and embed $G$ in $\R P^3$.
If either $G_1$ or $G_2$ contains all $0$-homologous cycles, then that subgraph
is crossing-change equivalent to a spatial embedding by Lemma \ref{cc}, and
hence
triple-linked by Theorem \ref{K7K7g}. Thus, in this case, $G$ has a
triple-link. Now suppose that both $G_1$ and $G_2$ contain a $1$-homologous
cycle. 

\vspace{\baselineskip}

Label the vertices of $G_1$ $\{1,2,3,4,5,6,7,A,B,C,D,E\}$ so that
$G[1,2,3,4,5,6,7]$ and $G[6,7,A,B,C,D,E]$ are isomorphic to $K_7$ and share
edge $(6,7)$. Up to isomorphism, there are three $3$-cycle equivalence classes
in $G_1$. The set $S = \{(1,2,3), (1,2,7)$, $(1,6,7)\}$ contains one
representative from each $3$-cycle equivalence class. By the same argument
given in Theorem \ref{K6K6RP3}, we may assume that $S$ contains a
$1$-homologous cycle, $c_1$. 

\vspace{\baselineskip}

If $G[A,B,C,D]$ contains a $1$-homologous cycle, then this cycle, $c_1$, and
the $1$-homologous cycle in $G_2$ form three disjoint $1$-homologous cycles, so
$G$ contains a triple-link. Otherwise, $G[A,B,C,D,E,6]$ contains a pair of
linked cycles by Lemma \ref{OHK4}. Following the proof in Theorem
\ref{K6K6RP3}, connect the linked cycle containing vertex $\{6\}$ to $c_1$ via
two disjoint paths. By Lemma \ref{BFLL}, $G$ contains a triple-link. \end{proof}

\vspace{\baselineskip}

The minor-minimality of the graph formed by identifying an edge of $K_7$ with
an edge from another copy of $K_7$ with respect to intrinsic triple-linking is
unknown in $\R^3$. If true, then the graph $G$ defined in Theorem \ref{K7e} is
also minor-minimal with respect to intrinsic triple-linking; a similar argument
to that in Theorem \ref{K6K6RP3} holds in this case as well. 

\vspace{\baselineskip}

We also remark that the graph $G(n)$ as defined in \cite{FFNP} is a
one-component minor-minimal intrinsically $(n+1)$-linked graph in $\R P^3$. The
arguments given in \cite{FFNP} hold in $\R P^3$ since $K_{4,4}-edge$ is
intrinsically linked in both $\R^3$ and $\R P^3$. 

\section{Graphs with linking number $\geq$ 1 in $\mathbb{R}P^3$}

In $\mathbb{R}P^3$, there are intrinsically linked graphs for which there exists
an embedding in which every pair of disjoint cycles has linking number less than
1, as a pair of linked cycles may have only one crossing.  Work has been done
\cite{F} in
$\mathbb{R}^3$ to find graphs containing disjoint cycles with large linking
number in every spatial embedding. Using the fact that $K_{10}$ is triple-linked
in $\mathbb{R}^3$, Flapan \cite{F} showed that every spatial embedding of
$K_{10}$ contains a two-component link $L\cup J$ such that, for some
orientation, $lk(L, J) \geq 2.$  A similar argument using Theorem 8 yields the
following proposition.

\begin{prop} {\it Every projective embedding of $K_{10}$ contains a
two-component link $L\cup J$ such that, for some orientation, $lk(L, J) \geq
1.$}
\end{prop}

It remains an open question to determine whether 10 is the smallest number for
which this property holds. At this point, we know the smallest $n$ is such that
$7 < n \leq 10$.

\vspace{\baselineskip}
\vspace{\baselineskip}

\end{document}